\documentclass[options]{article}
\usepackage{amsmath,amssymb}

\newcommand{\be}{\begin{equation}}
\newcommand{\ee}{\end{equation}}
\newcommand{\ba}{\begin{eqnarray}}
\newcommand{\ea}{\end{eqnarray}}
\newcommand{\baa}{\begin{eqnarray*}}
\newcommand{\eaa}{\end{eqnarray*}}
\newcommand{\bb}{\begin{thebibliography}}
\newcommand{\eb}{\end{thebibliography}}

\newcommand{\bi}[1]{\bibitem{#1}}
\newcommand{\lab}[1]{\label{#1}}
\newcommand{\re}[1]{(\ref{#1})}

\newcounter{my}
\newcommand{\he}%
   {\stepcounter{equation}\setcounter{my}%
   {\value{equation}}\setcounter{equation}0%
   }%
\newcommand{\she}%
   {\setcounter{equation}{\value{my}}%
    }%

\renewcommand\t{\tilde}

\newcommand\ve{\varepsilon}

\newtheorem{pr}{Proposition}

\newtheorem{lem}{Lemma}

\begin{document}

\vspace*{10mm}

\begin{center}

{\Large \bf Umbral "classical" polynomials}

\vspace{5mm}

{\large \bf Alexei Zhedanov}

\medskip

{\em Donetsk Institute for Physics and Technology, Donetsk 340114,  Ukraine}

\end{center}


\begin{abstract}
We study the umbral "classical"  orthogonal polynomials
with respect to a generalized derivative operator $\cal D$ which
acts on monomials as ${\cal D} x^n = \mu_n x^{n-1}$ with some
coefficients $\mu_n$. Let $P_n(x)$ be a set of orthogonal polynomials. Define the new polynomials
$Q_n(x) =\mu_{n+1}^{-1}{\cal D} P_{n+1}(x)$. We find necessary and sufficient conditions when the polynomials $Q_n(x)$ will also be orthogonal.  Apart from well known examples of the classical orthogonal polynomials we present a new example of umbral classical polynomials expressed in terms of elliptic functions.

\vspace{2cm}

{\it Keywords}: umbral classical polynomials; Hahn theorem

\vspace{2cm}

{\it AMS classification}: 05A40, 33C47, 42C05

\end{abstract}

\bigskip\bigskip


\newpage
\section{Introduction}
\setcounter{equation}{0} 
Let $P_n(x)=x^n + O(x^{n-1})$ be a set of monic orthogonal polynomials with respect to the linear functional $\sigma$
\be
\langle \sigma, P_n(x) P_m(x) \rangle = h_n \: \delta_{nm}, \quad n,m=0,1,2,\dots \lab{ort_P_sigma} \ee  
where $h_n$ are nonzero normalization constants.

The linear functional $\sigma$ acts on the space of polynomials. This action can be defined by the moments
\be
\langle \sigma, x^n \rangle = g_n, \lab{mom_def} \ee  
where $g_n, \: n=0,1,2,\dots$ is a sequence of complex numbers with the initial condition $g_0=1$. We assume that the linear functional $\sigma$ is nondegenerate. This means that all the Hankel determinants
\be
\Delta_n = \det|g_{i+k}|_{i,k=0}^{n-1} \lab{Hank_det} \ee
are nonzero $\Delta_n \ne 0, \; n=1,2, \dots$. There is the relation
\be
h_n = \frac{\Delta_{n+1}}{\Delta_n} \lab{h_n_H} \ee
from which it follows that the nondegenerate condition is equivalent to the condition $h_n \ne 0$.

Moreover, the orthogonal polynomials $P_n(x)$ satisfy the three-term recurrence relation \cite{Chi}
\be
P_{n+1}(x) + b_n P_n(x) + u_n P_{n-1}(x) = xP_n(x), \lab{3-term_P} \ee
where
\be
u_n = \frac{h_n}{h_{n-1}} \lab{u_h} \ee
Clearly, the nondegenerate condition is equivalent to the condition $u_n \ne 0, \: n=1,2,\dots$.

The monic polynomials 
\be
Q_n(x) = \frac{P_{n+1}'(x)}{n+1} = x^n + O(x^{n-1}) \lab{Q_P_drv} \ee
are in general, not orthogonal (here $f'(x)$ means the ordinary derivative). The set of polynomials $P_n(x)$ is called the {\it classical polynomials} if the polynomials $Q_n(x)$ are orthogonal polynomials with respect to another nondegenerate linear functional $\tau$
\be
\langle \tau, Q_n(x) Q_m(x) \rangle = \t h_n \: \delta_{nm}, \lab{ort_Q_sigma} \ee
where the linear functional $\tau$ is defined by the moments
\be
\t g_n = \langle \tau, x^n \rangle \lab{tg_def} \ee
and $\t h_n$ are new normalization constants.

The Hahn theorem \cite{Al-Salam} gives a complete classification of all classical polynomials. There are essentially 4 distinct types of such polynomials: Jacobi, Laguerre, Hermite and Bessel polynomials.

The ordinary derivative operator $\partial_x$ satisfies the property on monomials
\be
\partial_x x^n = n x^{n-1} \lab{drv_xn} \ee
This is defining property of the derivative operator. 

In the umbral calculus \cite{Roman} the formal derivative operator $\cal D$ can be introduced. Its action on monomials is defined by the formula
\be
{\cal D} x^n = \mu_n x^{n-1}, \quad n=0,1,2,\dots\lab{def_D} \ee
where $\mu_n$ is an arbitrary sequence of complex numbers with the restrictions $\mu_0=0$ and $\mu_n \ne 0, \: n=1,2,\dots$.  Clearly, the operator $\cal D$ decreases the degree of any polynomial by one. The obvious example of the formal derivative operator is the ordinary derivative operator  (in this case $\mu_n =n$). Another simple example is the q-derivative operator
\be
D_{q} f(x) =\frac{f(xq)-f(x)}{x(q-1)}. \lab{D_q} \ee
In this case $\mu_n = (q^n-1)/(q-1)$.

We say that the set of monic orthogonal polynomials $P_n(x), \: n=0,1,2,\dots $ satisfies the umbral classical property if the new set of monic polynomials
\be
Q_n(x) = \frac{{\cal D} P_{n+1}(x)}{\mu_{n+1}}, \quad n=0,1,2,\dots \lab{tPD} \ee
is another set of orthogonal polynomials (see e.g. \cite{Khol} for further details). We will adopt notation $\tau$ for the linear functional providing orthogonality of $Q_n(x)$ and $\t g_n$ for corresponding moments \re{tg_def}.

Hahn solved the problem of classification of the classical polynomials in cases of the operators $\partial_x$ and $D_q$ (see e.g. \cite{Al-Salam}).

The main purpose of the present paper is to find necessary and sufficient conditions for the polynomials $P_n(x)$ and for the sequence $\mu_n$ to satisfy the umbral classical property. We derive these conditions and analyze them presenting several special cases. Apart from well known examples of orthogonal polynomials belonging to the Hahn class (for both the ordinary and q-derivative operators) we present a new example of umbral claqssical polynomials connected with elliptic functions.

The paper is organized as follows.

In Section 2, we introduce the reciprocal operator $\cal R$ which transforms polynomials $Q_n(x)$ to $P_n(x)$. In general, this operator is not unique. However, if the polynomials $P_n(x)$ and $Q_n(x)$ are orthogonal then the operator $\cal R$ is completely determined. This allows to formulate a necessary and sufficient condition for polynomials $Q_n(x)$ to be orthogonal. From this condition it is possible to derive the eigenvalue problems $L P_n(x) = \lambda_n P_n(x)$ and $\t L Q_n(x) = \lambda_{n+1} Q_{n}(x)$ for polynomials $P_n(x)$ and $Q_n(x)$, where $L ={\cal R}{\cal D}, \;  \t L ={\cal D}{\cal R}$. The operator $\cal R$ acts on monomials as
$$
{\cal R}  x^n = \sum_{s=0}^{n+1} R_{ns} x^s
$$
with the condition $R_{n,n+1} \ne 0$.

In section 3, we study the simplest case when the matrix $R_{ns}$ has only 2 nonzero diagonals $R_{n,n+1}$ and $R_{nn}$, i.e. $R_{ns}=0$ when $s<n$. This leads to classical and q-classical polynomials (but not exhausts all classical polynomials)

In Section 4, we consider more general case when the matrix $R_{ns}$ has only finite number of diagonals, i.e. $R_{ns}=0$ when $s< n-j$ with fixed $j=1,2,\dots$. The main result of this section is that this condition is equivalent to the statement that $\mu_n$ satisfies a finite difference equation with constant coefficients.

In Section 5, we show that the polynomials $Q_n(x)$ and $P_n(x)$ are connected by a chain of spectral (i.e. Christoffel and Geronimus) transforms and then we obtain the reduced system of algebraic equations which is much simpler than inital one.  

In Section 6, we derive possible expression for the operators $\cal D$ and $\cal R$ in the case of local operator $\cal R$.  For the simplest case when $\cal R$ contains only 3 diagonals we show that the reduced system of equations leads to well known classical and q-classical polynomials.

Finally, in Section 7, we consider the case when the operator $\cal R$ is completely degenerate, i.e. $\lambda_n =1, \; n=1,2,\dots$. This leads to a nontrivial new example of umbral classical polynomials expressed in terms of elliptic functions. In a special case, when elliptic functions degenerate to rational functions, one obtains the Krall-Jacobi polynomials which appear to be umbral classical polynomials.

\section{Necessary and sufficient conditions} \setcounter{equation}{0} 
Let us introduce the new operator $\cal R$ which sends any polynomial of exact degree $n$ to a polynomial of exact degree $n+1$. Clearly, any such operator can be defined through it action on monomials
\be
{\cal R} x^n = \nu_{n+1} x^{n+1} + R_{nn} x^{n} + R_{n,n-1} x^{n-1} + \dots + R_{n0} \lab{R_x_n} \ee
with some coefficients $\nu_n, R_{nk}$. It is assumed that $\nu_n \ne 0, \; n=1,2,\dots$ (this is necessary for the condition $\deg({\cal R} x^n) = n+1$).

If orthogonal polynomials $P_n(x)$ are given by their expansion coefficients 
\be
P_n(x) = x^n + \sum_{s=0}^{n-1} A_{ns} x^s \lab{P_A_x} \ee
then the "derived" polynomials $Q_{n}(x)$ are defined uniquely
\be
Q_n(x) = x^n + \sum_{s=0}^{n-1} \frac{\mu_s A_{n+1,s}}{\mu_{n+1}} x^s. \lab{Q_A_x} \ee 
We can construct an operator $\cal R$ such that
\be
{\cal R} Q_n(x) = \nu_{n+1}  P_{n+1}(x). \lab{R_Q_P} \ee 
It is easily seen that there are infinitely many such operators. Indeed, the nonzero coefficients $\nu_1, \nu_2, \dots$ can be chosen arbitrarily while the coefficients $R_{nk}$ are defined uniquely from condition \re{R_Q_P}.

However, if the polynomials $Q_n(x)$ are orthogonal with respect to a nondegenerate functional $\tau$ then the operator $\cal R$ is defined uniquely, i.e. the coefficients $\nu_k, k=1,2,\dots$ are completely determined. Specifically, we have the
\begin{pr}
If the polynomials $Q_n(x)$ are orthogonal with respect to the nondegenerate linear functional $\tau$ then the relation
\be
\langle \tau, g(x) {\cal D} f(x) \rangle = \langle \sigma, f(x) {\cal R} g(x) \rangle \lab{gDf} \ee
holds for any pair of polynomials $f(x), g(x)$. Moreover, the coefficients $\nu_{n}$ are determined uniquely by the relation 
\be
\nu_{n+1} = \frac{\mu_{n+1} \t h_n }{h_{n+1}}, \quad n=0,1,\dots \lab{nu_h} \ee
\end{pr}  
{\it Proof}. It is sufficient to prove this proposition for $g(x) = Q_n(x), \; f(x) = P_{m+1}(x), n,m=0,1,\dots$. We have by orthogonality properties of the polynomials $Q_n(x)$ and $P_n(x)$ 
\be
\langle \tau, Q_n(x) {\cal D} P_{m+1} (x) \rangle = \mu_{n+1} \langle \tau, Q_n(x) Q_{m} (x) \rangle   =   \mu_{n+1} \t h_n \: \delta_{nm} \lab{QDP_nm} \ee
and
\be
\langle \sigma, P_{m+1}(x) {\cal R} Q_{n} (x) \rangle =  \nu_{n+1} \langle \sigma, P_{m+1}(x)  P_{n+1} (x) \rangle =  \nu_{n+1} h_{n+1} \: \delta_{nm} \lab{PRQ_nm} \ee
Hence the relation
\be
\langle \tau, Q_n(x) {\cal D} P_{m+1} (x) \rangle =\langle \sigma, P_{m+1}(x) {\cal R} Q_{n} (x) \rangle \lab{QDP_PRQ} \ee
holds for all $n,m=0,1,2,\dots$ provided that condition \re{nu_h} is fulfilled. By linearity, we obtain that relation \re{gDf} holds for any pair of polynomials $g(x),f(x)$ because we can expand $g(x)$ and $f(x)$ as finite linear combinations of polynomials $Q_n(x)$ and $P_n(x)$. This proves the proposition.

The inverse statement is also valid:
\begin{pr}
Assume there exists a linear nondegenerate functional $\tau$ such that relation \re{gDf} holds for any pair of polynomials $g(x), f(x)$. Then the polynomials $Q_n(x)$ are orthogonal with respect to the functional $\tau$ and the operator $\cal R$ is uniquely determined by this relation.
\end{pr}
In order to proof this proposition we again choose $g(x) = Q_n(x), \; f(x) = P_{m+1}(x), n,m=0,1,\dots$ and then \re{gDf} gives us the orthogonality property of the polynomials $Q_n(x)$
\be
\langle \tau, Q_n(x), Q_m(x) \rangle = \t h_n \: \delta_{nm} \lab{tau_ort_Q} \ee
We thus have that relation \re{gDf} is necessary and sufficient for polynomials $Q_n(x)$ to be orthogonal.

There is a simple but important consequence of these propositions.

Indeed, we have two relations between orthogonal polynomials $P_n(x)$ and $Q_n(x)$:
\begin{subequations}\label{PQ_rel}
\begin{align}
{\cal D} P_{n+1}(x) &= \mu_{n+1} Q_n(x) \lab{DPQ} \\
{\cal R} Q_n(x) &= \nu_{n+1} P_{n+1} (x) \lab{RQP1}
\end{align}
\end{subequations}
whence there are two eigenvalue problems
\be
L P_n(x) = \lambda_n P_n(x), \quad \tilde L Q_n(x) = \lambda_{n+1} Q_n(x), \lab{eig_PQ} \ee
where
\be
L= {\cal R} {\cal D}, \quad \tilde L = {\cal D} {\cal R} \lab{LL_def} \ee
and $\lambda_n=\mu_n \nu_n$.

From basic relation \re{gDf} it follows that the operator $L$ is symmetric with respect to the linear functional $\sigma$, i.e. for arbitrary polynomials $f(x),g(x)$ the relation
\be
\langle \sigma, f(x) L g(x) \rangle = \langle \sigma, g(x) L f(x) \rangle, \lab{sym_L} \ee
holds. Note that relation \re{sym_L} is necessary and sufficient in order for eignenpolynomials $P_n(x)$ to be orthogonal \cite{Zhe_hyp}. Similarly the operator $\t L$ is symmetric with respect to the functional $\tau$:
\be
\langle \tau, f(x) \t L g(x) \rangle = \langle \tau, g(x) \t L f(x) \rangle . \lab{sym_tL} \ee
 We thus have that "classical" property of the orthogonal polynomials $P_n(x)$ implies additional eigenvalues problems \re{eig_PQ}. In case when ${\cal D}$ is a differential operator with polynomial coefficients this property was established in \cite{KS} and developed in \cite{KY}. 

Additionally, one can assume the nondegenerate condition 
\be
\lambda_n \ne \lambda_m \lab{ndeg_lam} \ee
for all pairs of distinct $n$ and $m$. In this case the eigenvalue problem $L P_n(x) = \lambda_n P_n(x)$ defines all polynomials $P_n(x)$ uniquely (up to a common factor). However, in the last section we will consider the case with complete degeneration of the operator $L$.

In order to obtain an effective algebraic criterion of the orthogonality of the polynomials $Q_n(x)$  let us choose $g(x) = x^m, \: f(x) = x^n$ in \re{gDf}. Then we have the system of algebraic equations
\be
\mu_n \t g_{n+m-1} = \nu_{m+1} g_{n+m+1} + R_{mm} g_{n+m} + \dots R_{m0} g_n, \quad m,n=0,1,2,\dots \lab{main_mu_c} \ee 
But if relation \re{gDf} is valid for all monomials $g(x) = x^m, \: f(x) = x^n$, then this relation is valid for any pair of polynomials $f(x), g(x)$. 

We thus have the
\begin{pr}
System of algebraic equation \re{main_mu_c} is necessary and sufficient condition for polynomials $Q_n(x)$ to be orthogonal.
\end{pr}

These equations contains many unknowns: the moments $g_n$, the "derived" moments $\t g_n$, the coefficients $\mu_n$ and $\nu_n$ and the matrix elements $R_{nk}$ of the operator ${\cal R}$. One can exclude some of these unknowns.

Indeed, let us put $m=0$. Then \re{main_mu_c} becomes
\be
\mu_n \tilde g_{n-1} =  R_{00} g_n + \nu_{1} g_{n+1} \label{m_0_eq} \ee
and we can express the "derived" moments in terms of initial moments:
\be
\t g_n = \frac{R_{00} g_{n+1} + \nu_1 g_{n+2}}{\mu_{n+1}}, \quad n=0,1,2,\dots  \lab{tg_Rg} \ee
Expressions \re{tg_Rg} define $\t g_n$ uniquely because, by condition,  $\mu_n \ne 0, \; n=1,2,\dots$.

We can then present the remaining  equations \re{main_mu_c} in the form
\be
\mu_n (R_{00} g_{n+m} + \nu_1 g_{n+m+1}) = \mu_{n+m} \left(\nu_{m+1} g_{n+m+1} + R_{mm} g_{n+m} + \dots R_{m0} g_n \right), \quad m=1,2,\dots \lab{main_mod_1} \ee 
Moreover, putting $n=0$ in \re{main_mu_c} we get
\be
\nu_{m+1} g_{m+1} + R_{mm} g_{m} + \dots R_{m1} g_1 + R_{m0} =0 \lab{n=0_muc} \ee
Equations \re{n=0_muc} allow to determine the moments $g_1,g_2,g_3$ step-by-step in terms of the coefficients $R_{ns}$ of the operator $\cal R$. For example,
\be
g_1 = -\frac{R_{11}}{\nu_1}, \; g_2= -\frac{R_{11} g_1 + R_{10}}{\nu_2}, \dots  \lab{g_1_g_2} \ee 
The remaining unknowns are the coefficients $\mu_n, \nu_{n}, \; n=1,2,\dots$ and the matrix coefficients $R_{ns}$ of the operator $\cal R$.

However, general solution may appear to be rather complicated problem. Instead, we can assume some additional conditions on parameters $\mu_n$ or on the structure of the matrix $R$. This simplifies the problem and allows to find explicitly at least some of unknowns.

Before considering special cases, we note a simple 
\begin{pr} 
Let $\mu_n, \: g_n$ be a solution of the "umbral" classical polynpomials (initial conditions $\mu_0=0, g_0=1$ are asumed). 
Then $\alpha q^n \mu_n, \: p^n g_n$ with arbitrary nonzero parameters $\alpha, q, p$ is also a solution with the same initial conditions. 
\end{pr}
The proof of this proposition is almost obvious. Indeed, lets $S(\tau)$ be the dilation operator $S(\tau)P_n(x)= \tau^{-n} P_n(x\tau)$. Clearly, the polynomials $S(\tau) P_n(x)$ remain monic orthogonal polynomials. Applying this operator to the main condition  
\be
S(\tau_1) {\cal D} S^{-1}(\tau_2) S(\tau_2) P_n(x) = \mu_n S(\tau_1) Q_n(x) \lab{SSPQ} \ee
we achieve the statement of the proposition. 

Using this proposition we can combine all solutions into equivalence classes. For example a solution with $\mu_n = \kappa_1^n-\kappa_2^{n}$ with two different parameters $\kappa_1,\kappa_2$ is equvalent to the solution $\mu_n = 1-q^n$ with $q=\kappa_2/\kappa_1$.

\section{Two-diagonal operator $R$}
\setcounter{equation}{0}
In this section we consider a simplest special case when the matrix of the operator $R$ in the monomial basis has only two diagonals, i.e. we assume that the operator $\cal R$ acts as
\be
{\cal R} x^n = \nu_{n+1} x^{n+1} + \rho_n x^n \lab{2-diag_R_x} \ee
with some nonzero coefficients $\rho_n$.

Then we can present the system \re{main_mu_c} in the form
\be
\mu_n \tilde g_{n+m-1} = \nu_{m+1} g_{n+m+1} + \rho_m g_{n+m} \lab{2-diag_main} \ee
Putting $m=0$ we can eliminate $\tilde g_n$ from the system \re{2-diag_main}:
\be
\tilde g_n = \frac{\nu_1 g_{n+2}}{\mu_{n+1}} + \frac{\rho_0 g_{n+1}}{\mu_{n+1}} \lab{tc_2-diag} \ee
Then system \re{2-diag_main} reads
\be
\left( \frac{\nu_1 \mu_n}{\mu_{n+m}} - \nu_{m+1} \right) g_{n+m+1} + \left( \frac{\rho_0 \mu_n}{\mu_{n+m}} - \rho_{m} \right) g_{n+m} =0 \lab{red_2-diag-mom} \ee
For $m=0$ system \re{red_2-diag-mom} is trivial: $0=0$. The first nontrivial equations correspond to $m=1$:
\be
\left( \frac{\nu_1 \mu_n}{\mu_{n+1}} - \nu_{2} \right) g_{n+2} + \left( \frac{\rho_0 \mu_n}{\mu_{n+1}} - \rho_{1} \right) g_{n+1} =0 \lab{2-diag_m=1} \ee
and to $m=2$:
\be
\left( \frac{\nu_1 \mu_{n-1}}{\mu_{n+1}} - \nu_{3} \right) g_{n+2} + \left( \frac{\rho_0 \mu_{n-1}}{\mu_{n+1}} - \rho_{2} \right) g_{n+1} =0 \lab{2-diag_m=2} \ee
One can consider equations \re{2-diag_m=1}-\re{2-diag_m=2} as a linear homogeneous system with respect to unknowns $g_{n+1}, g_{n+2}$. Existence of nontrivial solution of such system is possible only if the determinant is equal to zero. This leads to the equation
\be
\alpha_1 \mu_{n} + \alpha_2 \mu_{n+1} + \alpha_3 \mu_{n+2} =0 \lab{lin_mu} \ee
Up to transformations $\mu_n \to \alpha q^n \mu_n, \: g_n \to p^n g_n$ there are only two distinct solutions of equation \re{lin_mu}:
\be
\mu_n = \frac{1-q^n}{1-q} \lab{mu_qn} \ee 
and
\be
\mu_n =n \lab{mu_n} \ee
The first case corresponds to the ordinary derivative operator, i.e. $\cal D= \partial_x$. The second case corresponds to the q-derivative operator $\cal D=D_q$.

It is easily verified that expressions for $\nu_n$ and $\rho_n$ are
\be
\nu_n = \alpha_1 + \alpha_2 q^{-n}, \quad \rho_n =  \alpha_3 + \alpha_4 q^{-n} \lab{nu_rho_q} \ee
in the case of solution \re{mu_qn}
and
\be
\nu_n = \alpha_1 + \alpha_2 n, \quad \rho_n =  \alpha_3 + \alpha_4 n \lab{nu_rho_n} \ee
in the case of solution \re{mu_n}. Here $\alpha_i, i=1,\dots,4$ are arbitrary parameters.

The operator $L = {\cal R}{\cal D}$ has the abstract "hypergeometric" form \cite{Zhe_hyp}
\be
L x^n = \lambda_n x^n + \tau_n x^{n-1} \lab{L_hyp} \ee
with 
\be
\lambda_n = \mu_n \nu_n, \quad \tau_n = \mu_n \nu_{n-1} \lab{lt_hyp} \ee
This corresponds to solutions obtained in \cite{Zhe_hyp}. In more details, the case \re{mu_qn} corresponds to the little q-Jacobi polynomials and their special and degenerate cases, while the case \re{mu_n} corresponds to the Jacobi polynomials and their degenerate cases - Laguerre and Bessel polynomials.

We thus see that the case of two-diagonal operator $\cal R$ corresponds to abstract "hypergeometric" polynomials classified in \cite{Zhe_hyp}. Note, nevertheless, that the class of solutions in \cite{Zhe_hyp} is wider: it includes so-called little -1 Jacobi polynomials. These polynomials do not appear in the case of two-diagonal operator $\cal R$ because for them 
\be
\mu_n = n + \eta(1-(-1)^n) \lab{Dunkl_mu} \ee
with some parameter $\eta$. This expression for $\mu_n$ does not appear in the case of two-diagonal operator $\cal R$.

\section{The local operator $R$}
\setcounter{equation}{0}
Results of the previous section can be generalized to the case when the matrix $R$ is finite-diagonal. This means that there exists a positive integer $j$ such that
\be
{\cal R} x^n = \nu_{n+1} x^{n+1} + K_n^{(0)} x^n + K_n^{(1)} x^{n-1} +  \dots + K_n^{(j-1)} x^{n-j+1}  \lab{j_Rx} \ee
where
\be
K_n^{(i)}= R_{n,n-i} \lab{K_R} \ee
There is an obvious truncation condition
\be
K_n^{(k)}=0, \quad k>n \lab{trunc_K} \ee
preventing appearing of terms with negative degrees in rhs of \re{j_Rx}

The corresponding operators $\cal R$ can be called the local operators, because they contain only a finite number of diagonals. Obvious examples of such operators are differential operators
\be
A_N(x) \partial_x^N +   A_{N-1}(x) \partial_x^{N-1} + \dots + A_0(x), \lab{R-diff} \ee  
where $N$ is a fixed positive integer and where $A_k(x)$ are polynomials in $x$ such that $\deg(A_k(x)) =k+1$

In the previous section we considered the simplest case of the local operator $j=1$ when the matrix $R$ is two-diagonal.

We have the following 
\begin{pr}
If the operator $\cal R$ has $j+1$ diagonals in the monomial basis \re{j_Rx} then the parameters $\mu_n$ satisfy the recurrence relation with constant coefficients 
\be
\alpha_0 \mu_n + \alpha_1 \mu_{n-1} + \dots + \alpha_{j+1} \mu_{n-j-1} =0, \quad n>j+1,  \lab{al_mu_rec} \ee
where the coefficients $\alpha_i, i=0,1,\dots,j+1$ do not depend on $n$.
\end{pr} 
{\it Proof}. From \re{j_Rx} and \re{main_mu_c} we have the system of $j+2$ homogeneous equations
\be
\mu_{n-i} \t g_{n+m-1} = \nu_{m+1+i} g_{n+m+1} + K_{m+i}^{(0)} g_{n+m} + K_{m+i}^{(1)} g_{n+m-1} + \dots + K_{m+i}^{(j-1)} g_{n+m-j+1}, \quad i=0,1,\dots,j+1 \lab{mom_i} \ee
for $j+2$ unknowns $\t g_{n+m}, g_{n+m+1}, g_{n+m}, \dots, g_{n+m-j+1}$.

This system should have a nonzero solution (otherwise the linear functionals are degenerate). Hence the determinant of this system should be equal to zero. This leads to the equation \re{al_mu_rec}. Note that the operator $\cal R$ with $j+1$ diagonals leads to the difference equation \re{al_mu_rec} of the same order $j+1$. For $j=1$ we obtain the case already considered in the previous section: two-diagonal operator $\cal R$ leads to second-order difference equation \re{lin_mu} for $\mu_n$.

There is an inverse statement with respect to above proposition:
\begin{pr}
Assume that $\mu_n$ satisfies linear difference equation \re{al_mu_rec} of the order $j+1$. Then the operator ${\cal R}$ should have no more than $j+2$ diagonals in the monomial basis:
\be
{\cal R} x^n = \nu_{n+1} x^{n+1} + K_n^{(0)} x^n + K_n^{(1)} x^{n-1} +  \dots + K_n^{(j)} x^{n-j} \lab{R_j+2} \ee
\end{pr} 
{\it Proof}. Assume that the coefficients $\mu_n$ satisfy relation \re{al_mu_rec}. We start with the basic conditions \re{main_mu_c}
\be
\mu_n \t g_{n+m-1} = \sum_{s=0}^{m+1} R_{ms} g_{n+s}, \lab{mu_g_R} \ee
where $R_{m,m+1}=\nu_{m+1}$. From \re{mu_g_R} we can obtain
\be
\sum_{i=0}^{j+1} \alpha_i \mu_{n-i} \t g_{n+m-1} = \sum_{i=0}^{j+1} \sum_{s=0}^{m+1} \alpha_i R_{m+i,s} g_{n+s-i} \lab{sum_alpha_g} \ee
By \re{al_mu_rec} the lhs of \re{sum_alpha_g} vanishes and we have the conditions
\be
\sum_{i=0}^{j+1} \sum_{s=0}^{m+1} \alpha_i R_{m+i,s} g_{n+s-i}=0 \lab{sum_alpha_g2} \ee
which should be valid for all $m=0,1,2,\dots$ and for all $n \ge m+1$. 

We can combine terms in front of moments:
\ba
&&\sum_{s=0}^{m+1} \sum_{i=0}^{j+1} \alpha_i R_{m+i,s} g_{n+s-i} = \nonumber \\ 
&&\alpha_{j+1} R_{m+j+1,0} g_{n-j-1}   + (\alpha_{j+1} R_{m+j+1,1} + \alpha_{j} R_{m+j,0}) g_{n-j} + \nonumber \\
&&(\alpha_{j+1} R_{m+j+1,2} + \alpha_{j} R_{m+j,1}+\alpha_{j-1} R_{m+j-1,0}) g_{n-j+1}\dots =0 \lab{cond_dep_mu} \ea 
The rhs of \re{cond_dep_mu} is a linear combination of the moments $g_{n}$ with coefficients not depending on $n$. From the nondegenerate condition it follows that all these coefficients should vanish:
\ba
&&\alpha_{j+1} R_{m+j+1,0} =  \alpha_{j+1} R_{m+j+1,1} + \alpha_{j} R_{m+j,0} = \nonumber \\
&&\alpha_{j+1} R_{m+j+1,2} + \alpha_{j} R_{m+j,1}\alpha_{j-1} R_{m+j-1,0} = \dots =0 \nonumber \ea 
whence
\be
R_{m0}=R_{m-1,1}= \dots =R_{m-i,i}=0, \quad i=0,1,\dots, j+1, \quad m \ge j+1 \lab{R_m=0} \ee
But condition \re{R_m=0} means that the matrix $R$ has only $j+2$ nonzero diagonals and we arrive at relation \re{j_Rx}.

Hence, when the coefficients $\mu_n$ satisfy linear difference equation \re{al_mu_rec} with constant coefficients, then basic equations \re{main_mu_c} can be presented in the form
\be
\mu_n \t g_{n+m-1} = \nu_{m+1} g_{n+m+1} + K_m^{(0)} g_{n+m} + \dots + K_m^{(j)} g_{n+m-j} \lab{main_c_j} \ee
Note the apparent "assymetry" between above Propositions. Indeed, if the operator $\cal R$ has $j+1$ diagonals then the coefficients $\mu_n$ satisfy the recurrence relation of the order $j+1$. However, if $\mu_n$ satisfies the recurrence relation of the order $j+1$ then the operator $\cal R$ has no more than $j+2$ diagonals. This can be explained by the observation that if the order of the recurrence relation for $\mu_n$ is $j+1$ then the operator $\cal R$ may have lesser diagonals than $j+2$, e.g. it may have $j+1$ diagonals.

\section{Reduced system of equations for the case of local operator $\cal R$}
\setcounter{equation}{0}
It is possible to eliminate terms $\t g_{n+m-1}$ from equations \re{main_c_j} by using the following 
\begin{lem}
Assume that the coefficients $\mu_n$ satisfy recurrence relation \re{al_mu_rec} of order $j+1$ and this order is minimal, i.e. any linear relation of the form 
$$
\sum_{i=0}^M \t \alpha_i \mu_{n-i}=0, \quad M<j+1
$$
is possible if and only if $\t \alpha_i=0, \: i=0,1,\dots ,M$.

Assume moreover, that the coefficients $\alpha_i$ in \re{al_mu_rec} satisfy the condition
\be
\sum_{i=0}^{j+1} \alpha_i  =0 \lab{sum_al_0} \ee
Then there exist constants $\beta_0, \beta_1, \dots, \beta_j$ such that 
\be
\sum_{i=0}^{j} \beta_i \mu_{n-i} =1 \lab{beta_mu_cond} \ee
\end{lem}
{\it Proof}. Assume that condition \re{sum_al_0} holds. Choose 
\be
\beta_0 = \gamma \: \alpha_0, \; \beta_i =\gamma \: \sum_{s=0}^{i} \alpha_i, \quad i=1,2,\dots, j \lab{beta_alpha} \ee 
with some complex constant $\gamma$.

Let us consider the expression 
\be
Y_{n} = \sum_{i=0}^j \beta_i \mu_{n-i} \lab{Y_mu} \ee 
It is easily verified that $Y_{n-1}=Y_n$ and hence $Y_n$ is a constant not depending on $n$. This constant cannot be zero, because otherwise $j+1$ is not the minimal order of equation \re{al_mu_rec}. Hence by an appropriate choice of the parameter $\gamma$ one can achieve condition \re{beta_mu_cond}.  This proves the Lemma.

The remaining question is: how to provide condition \re{sum_al_0} which is necessary for the Lemma. This can be achieved by the equivalence transformation $\mu_n \to q^{-n} \mu_n$ with some complex parameter $q$. Indeed, under this transformation we have $\alpha_i \to q^i \alpha_i$ and we should verify the condition
\be
\sum_{i=0}^{j+1} \alpha_i q^i =0 \lab{sum_al_q} \ee  
But \re{sum_al_q} is an algebraic equation of order $j+1$ with respect to the unknown $q$. It always has at least one nonzero solution. This means that by an appropriate equivalence transformation one can always achieve the desired condition \re{sum_al_0}.    

In what follows we will assume that condition \re{sum_al_q} holds. Then we have from \re{main_c_j}
\be
\sum_{i=0}^j \beta_i \mu_{n-i} \t g_{n+m-1} = \sum_{i=0}^j \beta_i \sum_{s=-1}^{j} K_{m+i}^{(s)} g_{n+m-s} \lab{sum_beta} \ee
where $K_m^{(-1)}=\nu_{m+1}$
Using \re{beta_mu_cond} we can present \re{sum_beta} in the form
\be
\t g_{n+m-1} = \sum_{s=-1}^j L_m^{(s)} g_{n+m-s}, \lab{tg_g} \ee
where
\be
L_m^{(s)} = \sum_{i=0}^j \beta_i K_{m+i}^{(s)} \lab{L_K} \ee
Putting $m=0$ in  \re{tg_g} we obtain that the modified moments $\t g_n$ are expressed as a linear combination of the  moments $g_n$
\be
\t g_n = \sum_{s=-1}^j L_0^{(s)} g_{n+1-s}, \quad n \ge j-1 \lab{tg_g0} \ee 
Putting $m=1$ we have similarly
\be
\t g_n = \sum_{s=-1}^j L_1^{(s)} g_{n+1-s}, \quad n \ge j-1 \lab{tg_g1} \ee 
Subtracting \re{tg_g0} and \re{tg_g1} we have
$$
\sum_{s=-1}^j (L_1^{(s)} - L_0^{(s)}) g_{n+1-s} =0
$$
Due to nondegenerate condition we obtain that $L_1^{(s)}=L_0^{(s)}$. Putting $m=2,3,\dots$ we have similarly that 
\be
L_0^{(s)}=L_1^{(s)}=\dots=L_j^{(s)} = \ve_s, \lab{L_ep} \ee
with some constants $\ve_s$, i.e. that the coefficients $L_i^{(s)}$ do not depend on $i$.

We thus have
\be
\t g_{n+j-1} = \sum_{s=-1}^j \ve_s g_{n+j-s}, \quad n=0,1,2,\dots \lab{tg_ep_g} \ee
This relation has a simple interpretation. Indeed, we can present \re{tg_ep_g} in the equivalent form
\be
x^{j-1} \tau = \pi_{j+1}(x) \sigma, \lab{Darb_t_s} \ee 
where 
\be
\pi_{j+1}(x) = \sum_{s=-1}^{j} \ve_s x^{j-s} \lab{pi_(x)} \ee
is a polynomial of degree $\le j+1$. The product $\pi(x) \sigma$ of the linear functional $\sigma$ by a polynomial $\pi(x)$ is defined as
\be
\langle \pi(x) \sigma, f(x) \rangle = \langle  \sigma, \pi(x)f(x) \rangle \lab{prod_pi} \ee
for any polynomial $f(x)$. On the other hand, the functional $\pi_{j+1}(x) \sigma$ corresponds to $j+1$ step Christoffel transform of the functional $\sigma$ \cite{ZheR}. Reciprocal to Christoffel is Geronimus transform \cite{ZheR}. We thus see that the functionals $\sigma$ and $\tau$ are related by sequences of Christoffel and Geronimus transforms (see \cite{Maroni}, \cite{ZheR} for further details).

Substituting \re{tg_ep_g} into \re{main_c_j} we obtain the reduced system of equations 
\be
\sum_{s=-1}^j (K_m^{(s)} - \mu_n \ve_s)g_{n+m-s} =0, \quad n,m=0,1,2,\dots  \lab{main_c_g} \ee
This system of algebraic equations is necessary and sufficient condition in order for polynomials $P_n(x)$ be umbral classical in case when the operator $\cal R$ is local.

We can further specify dependence of the coefficients $K_m^{(s)}, \: s=-1,0,\dots,j$ on $m$. Indeed, system of algebraic equations \re{main_c_g} is homogeneous in variables $g_{n+m+1}, g_{n+m},g_{n+m-1}, \dots, g_{n+m-j}$. Hence the determinant of the matrix (the rows are labelled by $i$, the columns are labelled by $s$)
\be
W(n,m)_{is}=K_{m+i}^{(s)}-\ve_s \mu_{n-i}, \quad  i=0,1,\dots,j+1, \; s=-1,0,\dots,j \lab{det_W} \ee 
should be zero
\be
\det(W(n,m))=0 \lab{W=0} \ee 
for all admissible values of $m,n$. Let us multiply each $i$-th row of the matrix $W(n,m)$ by $\alpha_i$ and add them to the last row. Then by \re{al_mu_rec}  the last row of the resulting matrix $\t W(n,m)$ will be
\be
\t W(n,m)\left|_{i=j+1} \right .   = \sum_{k=0}^j \alpha_k K_{m+k}^{(s)}, \; s=-1,0,1,\dots, j \lab{last_row} \ee
Hence all entries in the last row of the matrix $\t W(n,m)$ depend on $m$ only (but not on $n$) while entries of all other rows depend on both variables $n,m$. The determinant should vanish $\det(\t W(n,m))=0$ for all $n,m$. This is possible only if 
\be
\sum_{k=0}^j \alpha_k K_{m+k}^{(s)}=0, \quad s=-1,0,1,\dots,j \lab{rec_K} \ee 
Equations \re{rec_K} have simple meaning: all the coefficients $K_m^{(-1)}=\nu_{m+1}, \: K_m^{(0)}, \dots, K_m^{(j)}$ are solutions of difference equation \re{rec_K} with constant coefficients $\alpha_k$. This equation has the same order $j+1$ as equation \re{main_c_g} for $\mu_n$ but is "reverse" with respect to \re{main_c_g}. For example, for $j=1$ we have equation for $\mu_n$
\be
\alpha_0 \mu_n + \alpha_1 \mu_{n-1} + \alpha_2 \mu_{n-2}=0
\lab{eq_mu_j=1} \ee 
and corresponding equation
\be
\alpha_2 K_{n}^{(s)} + \alpha_1 K_{n-1}^{(s)} + \alpha_0 K_{n-2}^{(s)} =0, \quad s=-1,0,1 \lab{eq_K_j=1} \ee
for $K_n^{(s)}$.

This means that in the case of local operators $\cal R$ the matrix coefficients $\mu_n, \: K_{n}^{(s)}$ of the operators $\cal D$ and $\cal R$ are solutions of difference equations of the same order with constant coefficients.

Return now to the reduced system \re{main_c_g} with fixed $j$. This system contains moments $g_n$ as main unknowns. The coefficients $\mu_n$ and $K_n^{(s)}$ can be explicitly presented as solutions of difference equations \re{al_mu_rec} and \re{rec_K}.

Solutions of such equations are well known (see e.g. \cite{Lando}). In the non-degenerate case the generic solution is
\be
\mu_n = a_1 q_1^n + a_2 q_2^n + \dots a_{j+1} q_{j+1}^n, \lab{gen_sol_mu} \ee 
where $q_i, \: i=1,2\dots,j+1$ are distinct roots of the Euler characteristic equation
\be
\alpha_0 q^{j+1} + \alpha_1 q^j + \dots + \alpha_{j} q +\alpha_{j+1}=0 \lab{char_q} \ee
and where $a_k, \: k=1,2,\dots,j+1$ are arbitrary parameters with the restriction
$$
a_1+a_2 + \dots + a_{j+1}=0
$$
in order to satisfy the initial condition $\mu_0=0$.

Similarly, for the coefficients $K_{n}^{(s)}$ we have the expressions
\be
K_n^{(s)} = r_1^{(s)} q_1^{-n} + r_2^{(s)} q_2^{-n} + \dots + r_{j+1}^{(s)} q_{j+1}^{-n} \lab{K_n_expl} \ee
with some constants $r_k^{(s)}$ and with the same characteristic roots $q_k$. Notice that in \re{K_n_expl} we have linear combination of negative degrees $q_k^{-n}$ in contrast to expression \re{gen_sol_mu}. The reason is that characteristic equation for the coefficients $K_n^{(s)}$ is 
\be
\alpha_{j+1} q^{j+1} + \alpha_j q^j + \dots + \alpha_{1} q +\alpha_{0}=0 \lab{char_q_K} \ee 
and it is clear that equation \re{char_q_K} has the roots $q_1^{-1}, q_2^{-1}, \dots, q_{j+1}^{-1}$ in case if equation \re{char_q} has the roots $q_1, q_2, \dots, q_{j+1}$.

Substituting these expressions into \re{main_c_g} we obtain a system of $j+2$ equations for the unknown moments $g_n$. Each equation of this system is a difference equation of order $j+1$ with respect to the same unknown sequence $g_n$. Compatibility conditions for these equations should give restrictions upon the parameters $q_i, a_i, r_i^{(s)}, \ve_s$. For $j=1$ these restrictions are very simple and will be analyzed in the next section. However for $j>1$ the detailed analysis of these restrictions could be a nontrivial problem.

When some of the roots $q_i$ coincide with one another solution for $\mu_n$ can contain terms like $n^k, k=0,1,\dots$ together with terms $q^k$. Analysis of this (degenerate) can be done in the same manner.

\section{Explicit expressions for the operators $\cal D$ and $\cal R$}
\setcounter{equation}{0}
What about possible expression of the operators $\cal D$ and $\cal R$ in the case of local operator $\cal R$? 

Using explicit expression \re{gen_sol_mu} for $\mu_n$ in the nondegenerate case we see that the operator $\cal D$ can be presented in the form 
\be
{\cal D} = x^{-1} \left(a_1 T_1 + a_2 T_2 + \dots a_{j+1} T_{j+1} \right), \lab{D_finite} \ee
where
\be
T_k f(x) = f(q_k x), \quad k=1,2,\dots, j+1 \lab{T_k} \ee
is the dilation operator. As a special case (when $j=1$) we have the q-derivative operator $D_q$ (in this case $q_1=q, q_2=1$).

Similarly, from solution \re{K_n_expl} we arrive at expression of the operator $\cal R$
\be
{\cal R} = x^{-j} \left( V_{j+1}^{(1)}(x)T_1^{-1}  + V_{j+1}^{(2)}(x)T_2^{-1} + \dots + V_{j+1}^{(j+1)}(x)T_{j+1}^{-1}\right), \lab{R_T} \ee
where $V_{j+1}^{(k)}(x)$ are polynomials of degree $j+1$ and where the operator $T_k^{-1}$ is defined as 
\be
T_k^{-1} f(x) = f(q_k^{-1} x). \lab{T_k-} \ee

In the degenerate case, when some of the rots $q_k$ coincide with one another, we obtain a combination of q-difference and ordinary derivative operators. In a special case, when all $q_k=1$ we obtain the differential operator  
\be
{\cal D} = \sum_{k=0}^{j-1} {\gamma_k x^k \partial_x^{k+1}} \lab{D_differential} \ee
with some constants $\gamma_k$.

A more interesting situation occurs when some of the characteristic roots degenerate to $q_k=-1$. In this case one can expect appearing of the Dunkl type operators which contain the reflection operator $R$ defined as $Rf(x)=f(-x)$.

Consider e.g. the difference equation \re{al_mu_rec} for $\mu_n$ of the 3-rd order (i.e. $j=2$)
\be
\mu_n - \mu_{n-1} -\mu_{n-2} + \mu_{n-3} =0 \lab{rec_mu_Dunkl} \ee
Generic solution of \re{rec_mu_Dunkl} with the initial condition $\mu_0=0$ is
\be
\mu_n = \kappa \left( n + \eta (1-(-1)^n) \right), \lab{Dunkl_sol_mu} \ee
where $\kappa, \eta$ are arbitrary parameters. The parameter $\kappa$ is the common factor and it is sufficient to put $\kappa=1$. We then obtain expression \re{Dunkl_mu} which corresponds to the Dunkl operator
\be
{\cal D} = \partial_x +  \eta x^{-1} (1-R). \lab{Dunkl_D} \ee
In \cite{VZ_missing} it was shown that the little -1 Jacobi polynomials are umbral classical with respect to the Dunkl operator \re{Dunkl_D}.

We thus see that the case of the local operators $\cal R$ leads to the operators $\cal D$ which can be considered as a natural generalization of the derivative  or q-derivative operators.

Consider the simplest example when $j=1$. We already know that in this case (to within equivalence transforms) there are two essentially different solutions for $\mu_n$: $\mu_n=n$ and $\mu_n = 1-q^n$. In both cases condition \re{sum_al_0} holds, hence we can apply formulas \re{Darb_t_s} and \re{main_c_g}. 

Formula \re{Darb_t_s} reads
\be
\tau = \pi_2(x) \sigma \lab{tau_Chr_sigma} \ee
Formula \re{tau_Chr_sigma} means that the "derived" functional $\tau$ is obtained from the initial functional $\sigma$ by application of at most two Christoffel transforms. 

Formula \re{main_c_g} now reads
\be
(\nu_{m+1}-\ve_{-1} \mu_n)g_{n+m+1} + (K_m^{(0)}-\ve_{0} \mu_n)g_{n+m} + (K_m^{(1)}-\ve_{1} \mu_n)g_{n+m-1}=0 \lab{j=1_main} \ee
Consider first the case $\mu_n=n$. This corresponds to the derivative operator ${\cal D}=\partial_x$. The difference equation of minimal order for $\mu_n$ is
\be
\mu_n - 2 \mu_{n-1} + \mu_{n-2}=0 \lab{mu_eq_class} \ee
Equation for the coefficients $K_n^{(s)},\: s=-1,0,1$ is the same
\be
K_n^{(s)}  - 2 K_{n-1}^{(s)} + K_{n-2}^{(s)} =0 \lab{K_eq_clas} \ee
whence generic solution is 
\be
K_n^{(s)}= \xi_s n + \eta_s, \quad s=-1,0,1 \lab{K_n_class} \ee
with arbitrary constants $\xi_s, \eta_s$.
From the boundary condition $K_0^{(1)}=0$ we have $\eta_1=0$.

We thus have the equation
\be
(\xi_{-1} m + \eta_{-1} -\ve_{-1} n)g_{n+m+1} + (\xi_{0} m + \eta_{0} -\ve_{-1} n)g_{n+m} + (\xi_{1} m  -\ve_{1} n)g_{n+m-1}=0 \lab{eq_g_clas} \ee
This equation should be valid for all admissible pairs $m,n=0,1,2,\dots$. Clearly, this is possible only under the restrictions
\be
\ve_s=-\xi_s, \quad s=-1,0,1 \lab{res_ep_xi_clas} \ee 
Then equations \re{eq_g_clas} become
\be
(\xi_{-1} n + \eta_{-1} )g_{n+1} + (\xi_{0} n + \eta_{0} )g_{n} +  \xi_{1} n g_{n-1}=0 \lab{eq_g_clas_2} \ee
with arbitrary parameters $\xi_s, \eta_s$. Equation \re{eq_g_clas_2} was first derived by Geronimus \cite{Ger1} in his studying of the Hahn problem. This equation determines all classical orthogonal polynomials: Jacobi, Laguerre, Hermite and Bessel.   

Another way to see this is to consider the operator $L={\cal R} {\cal D}$. In our case ${\cal D}=\partial_x$ and from \re{K_n_class} we can reconstruct the operator {\cal R} explicitly 
\be
{\cal R} = (\xi_{-1} x^2 + \xi_0 x+\xi_1) \partial_x + \eta_{-1}x + \eta_0. \lab{R_hyp_clas} \ee
The operator $L$ coincides with generic hypergeometric operator 
\be
L= (\xi_{-1} x^2 + \xi_0 x+\xi_1) \partial_x^2 + (\eta_{-1}x + \eta_0)\partial_x .\lab{L_hyp_clas} \ee
The corresponding eigenvalue problem 
\be
L P_n(x) = \lambda_n P_n(x) \lab{LP_clas} \ee
is known to generate all classical orthogonal polynomials \cite{Al-Salam}, \cite{NSU}.

Consider the case $\mu_n= 1-q^n$. The minimal difference equation \re{al_mu_rec} with the condition $\alpha_0+\alpha_1+\alpha_2=0$ is
\be
\mu_{n} -(q+1)\mu_{n-1} + q \mu_{n-2}=0 \lab{min_eq_mu_q} \ee
Corresponding equation for $K_n^{(s)}, \: s=-1,0,1$ is obtained from \re{min_eq_mu_q} by the change $q \to q^{-1}$:
\be
K_{n}^{(s)} -(1+q^{-1})K_{n-1}^{(s)} + q^{-1} K_{n-2}^{(s)} =0 \lab{min_eq_K_q} \ee
Whence
\be
K_{n}^{(s)}= \xi_s q^{-n} + \eta_s, \quad \lab{K_q_class} \ee
From the boundary condition $K_0^{(1)}=0$ it follows that $\eta_1=-\xi_1$.

Equation \re{j=1_main} becomes
\be
(\xi_{-1} q^{-m} +\eta_{-1} -\ve_{-1} (1-q^n))g_{n+m+1} + (\xi_{0} q^{-m} +\eta_{0} -\ve_{0} (1-q^n))g_{n+m} + (\xi_{1} q^{-m} -\xi_{1} -\ve_{1} (1-q^n))g_{n+m-1}=0 \lab{eq_g_q} \ee
Equations \re{eq_g_q} should be compatible for all admissible $n$ and $m$. This leads to the restrictions 
\be
\ve_s=\eta_s, \quad s=-1,0,1. \lab{res_e_eta_q} \ee
Then we arrive at equation for the moments
\be
(\xi_{-1} +\eta_{-1} q^n) g_{n+1} + (\xi_{0} +\eta_{0} q^n) g_{n} + \xi_{1}(1-q^n) g_{n-1} =0, \quad, n=0,1,2,\dots \lab{mom_eq_q} \ee
with arbitrary parameters $\xi_s, \eta_s$. This equation determines uniquely all moments $g_1,g_2, \dots$ starting with $g_0=1$. 

In order to find corresponding orthogonal polynomials $P_n(x)$ let us consider the operator $\cal R$. From \re{R_j+2} and \re{K_q_class} we have
\be
{\cal R} = (\xi_{-1} x + \xi_0 + \xi_1 x^{-1})T_q^{-1} + \eta_{1} x  + \eta_0 - \xi_1 x^{-1}, \lab{R_q_clas} \ee
where $T_q^{-1}f(x) = f(x/q)$.

The operator $\cal D$ is
\be
{\cal D} = x^{-1}(1-T_q). \lab{D_q_clas} \ee
It is easily verified that the operator $L = {\cal R}{\cal D}$ coincides with the operator which generates the difference equation $LP_n(x) = \lambda_n P_n(x)$ for the big q-Jacobi polynomials  \cite{KLS}, \cite{Ismail}. For special choices of the parameters $\xi_i, \eta_i$ we obtain all "q-classical" polynomials on the q-exponential grid $x_s = q^s$ \cite{NSU} found by Hahn \cite{Hahn}.  

We thus demonstrated that the case $j=1$ of the local operator $\cal R$ leads to all classical and q-classical polynomials. Corresponding operators $\cal D$ are either ordinary derivative $\partial_x$ or q-derivative $D_q$ operators.

An interesting open problem is finding nontrivial examples of umbral classical polynomials for $j>1$. For example, one can try to find all Dunkl classical polynomials with respect to the Dunkl operator \re{Dunkl_D}. Such problem was completely solved  in \cite{Cheikh} for the special case of symmetric polynomials, i.e. when $P_n(-x) = (-1)^n P_n(x)$. In \cite{VZ_missing} it was shown that the non-symmetric little -1 Jacobi polynomials are Dunkl classical. However, the problem of describing all Dunkl classical polynomials is still open.

\section{The degenerate case}
\setcounter{equation}{0}
There is another case when the problem of "umbral" classical polynomials can be solved explicitly. This corresponds to the full degeneration of the eigenvalue problems \re{eig_PQ} when $\lambda_0=0$ and $\lambda_n =1$ for $n=1,2,\dots$. In this case the operator  $\tilde L$ is reduced to the identity operator on polynomials $Q_n(x)$ 
\be 
\tilde L Q_n(x) = Q_n(x), \quad n=0,1,2,\dots  \lab{tL_id_Q} \ee 
But the polynomials $Q_n(x)$ form a basis on the linear space of polynomials. Hence the operator $\tilde L$ is the identity operator on all polynomials:
\be
\tilde L x^n = x^n \lab{tL_id} \ee    
From condition \re{tL_id} one can derive that
\be R_{ni}=0, \quad i=1,2,3,\dots, n, \quad R_{n,n+1}=\nu_{n+1}=\frac{1}{\mu_{n+1}} \lab{R=0} \ee
This means that the operator $\cal R$ acts on monomials as
\be
{\cal R} x^n = \frac{x^{n+1}}{\mu_{n+1}} + R_{n0}, \quad n=0,1,2,\dots \lab{deg_Rx} \ee  
with some unknown coefficients $R_{n0}$. In order to relate these coefficients with moments  $g_n$ of the functional $\sigma$ let us consider the operator $L = {\cal R} {\cal D}$. By \re{eig_PQ} we have
\be
L P_n(x) = \left\{ {0, \quad n=0 \atop P_n(x), \quad n=1,2,3,\dots}    \right . \lab{LP_deg} \ee
From the orthogonality relation $\langle \sigma, P_n(x) \rangle =0, \; n=1,2,\dots$ it follows that
\be
\langle \sigma, L f(x) \rangle =0 \lab{Lf0} \ee
for any polynomial $f(x)$. In particular,
\be
\langle \sigma, L x^n \rangle =0, \quad n=0,1,2,\dots \lab{Lx0} \ee
But the operator $L$ acts on monomials as
\be
L x^n = \mu_n {\cal R} x^{n-1} = \left\{ 
  \begin{array}{l l}
    0 & \quad \text{if $n=0$}\\
    x^n + \mu_n R_{n-1,0} & \quad \text{if $n \ge 1$}
  \end{array} \right. \lab{Lxn_deg} \ee 
From \re{Lx0} we find that
\be
R_{n0}= -\frac{g_{n+1}}{\mu_{n+1}}, \quad n=0,1,2,\dots \lab{Rn0} \ee 
Hence the action of the operator $\cal R$ on monomials is
\be
{\cal R} x^n = \frac{1}{\mu_{n+1}} x^{n+1} - \frac{g_{n+1}}{\mu_{n+1}} \lab{R_x^n_3} \ee
Thus  relations \re{main_mu_c} for the degenerate case can be presented in the form
\be
\mu_n \mu_{m+1} \tilde g_{n+m-1} = g_{n+m+1} - g_{m+1} g_n, \quad n,m=0,1,2,\dots \lab{red_deg_mu} \ee
Putting $m=0$ in \re{red_deg_mu} we can eliminate $\tilde g_n$:
\be
\tilde g_n = \frac{g_{n+2} - g_1 g_{n+1}}{\mu_1 \mu_{n+1}} \lab{tg_c} \ee
Then \re{red_deg_mu} reads
\be
\frac{\mu_n \mu_{m}}{\mu_1 \mu_{n+m-1}}   = \frac{g_{n+m} - g_{m} g_n}{g_{n+m} - g_1 g_{n+m-1}}, \quad m,n=1,2,3,
\dots  \lab{red_mu_c} \ee
Note an obvious symmetry of \re{red_mu_c}  with respect to the permutation $m \leftrightarrow n$. As expected, equations \re{red_mu_c} are invariant with respect to the equivalence transformations
\be
\mu_n \to \alpha q^n \mu_n, \quad g_n \to p^n g_n \lab{tr_mu_c} \ee
with arbitrary nonzero parameters $\alpha, q, p$.

In order to find solutions of equations \re{red_mu_c} we observe that if $\mu_n$ and $g_n$ are analytic functions of their arguments, i.e. $\mu_n = \mu(n), g_n = g(n)$ then \re{red_mu_c} becomes the functional equation
\be
g(n+m) = g(n)g(m) + \mu(n) \mu(m) \psi(n+m), \lab{BC_fe} \ee
where
\be
\psi(n) = \frac{g(n)-g(1)g(n-1)}{\mu(1) \mu(n-1)} \lab{psi_def} \ee
Remarkably, equation \re{BC_fe} belongs to a class of functional equations whose general solution can be expressed in terms of the Weierstrass sigma functions \cite{BC}, \cite{BB}. Up to transformations \re{tr_mu_c} we can present this solution as 
\be
\mu_n = \frac{\sigma(wn)}{\sigma(w(n+\alpha))}, \quad g_n = \frac{\sigma(w \alpha)}{\sigma(w \beta)}\frac{\sigma(w(n+\beta))}{\sigma(w(n+\alpha))}, \lab{gen_el_sol} \ee
where $w,\alpha,\beta$ is an arbitrary nonzero parameters ($\beta \ne \alpha$) and the Weierstrass sigma function $\sigma(z) \equiv \sigma(z; {\mathfrak g_2}, {\mathfrak g_3})$ is defined by the standard way \cite{WW}. The elliptic parameters ${\mathfrak g_2},{\mathfrak g_3}$ can be arbitrary; we omit them in notation $\sigma(z)$ for brevity.

For the modified moments one obtains
\be
\t g_n = \frac{\sigma(w \alpha) \sigma(w(\beta-\alpha))  \sigma(w(n+\alpha+\beta+2))}{\sigma^2(w \beta)\sigma(w(n+\alpha+2))} \lab{tg_ell} \ee

One can check directly that \re{gen_el_sol} and \re{tg_ell} is a solution of equations \re{red_mu_c} and \re{tg_c} using well known Riemann identity for sigma function  \cite{WW}. It is clear from \re{gen_el_sol} that the initial conditions $\mu_0=0, \: g_0 =1$ are valid. However, $\t g_0 \ne 1$. This is not essential because all moments are defined up to a common factor: transformation $\t g_n \to \kappa \t g_n$ leads to the same orthogonal polynomials $Q_n(x)$.

From explicit expression \re{gen_el_sol} for the moment it is possible to identify the orthogonal polynomials $P_n(x)$ coincide with those considered in \cite{TZ_dense}. The latter can be obtained from elliptic solutions of the qd-algorithm. They have the expression
\be
P_n(x) = B_n \: {_3}E_2 \left( {-n, \alpha+n, 1+\alpha-\beta -n(\alpha+n) \atop \alpha, \alpha-\beta -n(\alpha+n) }; x\right) \lab{P_ell} \ee
Here ${_3}E_2(x)$ stands for the elliptic hypergeometric function. The coefficient $B_n$ are chosen to satisfy the property $P_n(x)=x^n + O(x^{n-1})$. 

The elliptic hypergeometric function is defined as \cite{TZ_dense}
\be
{_3}E_2\left({-n,a_1,a_2 \atop b_1, b_2 };x\right) =\sum_{k=0}^n \frac{[-n]_k [a_1]_k [a_2]_k}{[1]_k [b_1]_k [b_2]_k} x^k \lab{E_hyp} \ee 
and
\be
[a]_k = y(a)y(a+1) \dots y(a+k-1) \lab{el_Pochh} \ee
is the elliptic Pochhammer symbol.  Here we denote $y(x) \equiv \sigma(wx)$ for brevity. Note that our definition of the elliptic hypergeometric function is slightly different from the conventional one (see, e.g. \cite{Spi}).

The polynomials $P_n(x)$ satisfy the three-term recurrence relation
\be
P_{n+1}(x) -(A_n+C_n)P_n(x) + A_{n-1}C_n P_{n-1}(x) =xP_n(x), \quad P_0=1, P_{-1}=0, \lab{3term_ell} \ee 
where
\be
A_n = \frac{y^2(n+\alpha) y(\beta+\alpha n + (n+1)^2) y(\beta+ \alpha (n-1) + n(n-1))}{y(2n+\alpha) y(2n+\alpha+1) y(\beta+\alpha (n-1) + n^2) y(\beta+ \alpha n + n(n+1))} \lab{A_n_ell} \ee
and
\be
C_n = \frac{y^2(n) y(\beta+\alpha (n-2) + (n-1)^2) y(\beta+ \alpha n + n(n+1))}{y(2n+\alpha) y(2n+\alpha-1) y(\beta+\alpha (n-1) + n^2) y(\beta+ \alpha (n-1) + n(n-1))} \lab{C_n_ell} \ee
It is seen that for all possible values of the parameters $\alpha,\beta,w$ the coefficient $u_n$ cannot be positive for all $n$. This means that there is no positive orthogonality measure on the real line for the polynomials $P_n(x)$. Nevertheless, the polynomials $P_n(x)$ are nondegenerate for generic choice of these parameters.

The derived polynomials $Q_n(x)$ differ from the polynomials $P_n(x)$ by the change of the parameters $\alpha \to \alpha+2, \; \beta \to \beta+\alpha+2$. This means that the polynomials $Q_n(x)$ again satisfy the umbral classical property, and this process can be continued infinitely.

When ${\mathfrak g_2}={\mathfrak g_3}=0$ the Weiersrtrass elliptic functions degenerate to the rational ones, in particular, $\sigma(x)=x$. In this case the moments become
\be
g_n = \frac{\alpha}{\beta} \: \frac{n+\beta}{n+\alpha} = \frac{\alpha}{\beta}\left( 1+ \frac{\beta-\alpha}{n+\alpha} \right) \lab{g_n_clas} \ee   
It is clear from \re{g_n_clas} that the moments admit the representation
\be
g_n = \frac{\alpha(\beta-\alpha)}{\beta} \int_0^1 x^n w(x) dx , \lab{g_n_w} \ee
where
\be
w(x) = x^{\alpha-1} + \frac{1}{\beta-\alpha} \delta(x-1) \lab{w_Krall} \ee
and $\delta(x)$ is the Dirac delta-function. From \re{w_Krall} it is seen that the orthogonality measure contains absolute continuous part $x^{\alpha-1}$ which corresponds to special Jacobi polynomials $P_n^{(\alpha-1,0)(x)}$ (orthognal on the interval $[0,1]$) and a concentrated mass at point $x=1$. Orthogonal polynomials corresponding to the measure \re{w_Krall} are known as the Krall-Jacobi polynomials \cite{Little}. They are remarkable because of their bispectrality property: they are eigenfunctons of a linear fourth-order differential operator \cite{Little2}. In this rational limit the elliptic hypergeometric function in \re{P_ell} becomes the ordinary hypergeometric function ${_3}F_2(x)$.

The recurrence relation for the Krall-Jacobi polynomials is \re{3term_ell}, in expressions  \re{A_n_ell} and \re{C_n_ell} one should put $y(n)=n$ which corresponds to the classical limit $\sigma(x) \to x$ of the elliptic functions.

Thus the case of complete degeneration of the operators $L$ and $\t L$ leads to nontrivial examples of orthogonal polynomials $P_n(x)$ with the umbral classical property \re{tPD}. Note that from \re{deg_Rx} it follows that  the operator $\cal R$ is non-local in this case. In some sense this is the simplest example of the nonlocal operator $\cal R$: apart from the term $\nu_{n+1} x^{n+1}$ there is only a constant term $R_{n0}$.  An interesting open problem is to find other nontrivial examples of umbral classical polynomials with nonlocal operators $\cal R$.

\section{Conclusions}
\setcounter{equation}{0}
The main result of the paper is the system of algebraic relations  \re{main_mu_c} which is necessary and sufficient condition for umbral classical polynomials $P_n(x)$. Solution of this system (i.e. finding explicit expressions for $\mu_n, \nu_n, g_n$ etc) is rather nontrivial problem. A natural restriction upon the structure of the operator $\cal R$ - namely, the local property - leads to much simpler reduced system of conditions  \re{main_c_g}.   In turn, this local property is equivalent   to the simple difference equation \re{al_mu_rec} with constant coefficients for the unknowns $\mu_n$. This leads to the explicit form \re{D_finite} of the operator $\cal D$ being a generalization of the q-derivative operator. In degenerate cases one can obtain pure differential or Dunkl type operators. 

Moreover, we have shown that for the case of local operators $\cal R$ the polynomials $P_n(x)$ and $Q_n(x)$ are related by a chain of Darboux transforms \re{Darb_t_s}.

The solution depends on the number $j+2$ of the diagonals of the matrix corresponding to the operator $\cal R$. Equivalently, $j+1$ is the order of the finite difference equation \re{main_mu_c}. The simplest case $j=1$ leads to either classical (i.e. Jacobi, Laguerre, Hermite and Bessel) or q-classical polynomials on the q-exponential grid. Already the next case $j=2$ is far from being studied in details.  We already know that the Dunkl classical polynomials belong to the class $j=2$. Classification of all orthogonal polynomials which possess the umbral classical property for $j=2$ would be an interesting open problem.

Moreover, we have shown that the local operators $\cal R$ do not exhaust all possible cases of umbral classical polynomials. There is a nontrivial example of the nonlocal operator $\cal R$ leading to polynomials expressed in terms of elliptic functions. One can expect existing of the more general solutions for the nonlocal operator $\cal R$.

\bigskip\bigskip
{\Large\bf Acknowledgments}
\bigskip

\noindent The author is grateful to V. Spiridonov, S.Tsujimoto and L.Vinet for
discussion.

\newpage

\bb{99}


\bi{Al-Salam} W.A. Al-Salam, {\it Characterization theorems for orthogonal polynomials}, in: P. Nevai (Ed.), Orthogonal Polynomials: Theory and Practice, NATO ASI Series C: Mathematical and Physical Sciences, vol. {\bf 294}, Kluwer Academic Publishers, Dordrecht, pp. 1–-24.

\bi{Cheikh} Y. Ben Cheikh and M.Gaied, {\it Characterization of the Dunkl-classical symmetric orthogonal polynomials}, Appl. Math. and Comput.
{\bf 187}, (2007) 105--114.




\bi{BB} H.W.Braden and V.M.Buchstaber, {\it The general analytic solution of a functional equation of addition type}, SIAM J. Math. Anal., {\bf 28} (1997), 903--923.

\bi{BC} M. Bruschi and F. Calogero, {\it General analytic solution of certain functional equations of
addition type}, SIAM J. Math. Anal., {\bf 21} (1990), 1019--1030.


\bi{Chi} T. Chihara, {\it An Introduction to Orthogonal
Polynomials}, Gordon and Breach, NY, 1978.


\bi{Ger1} Ya.L.Geronimus, {\it On polynomials orthogonal with respect to to
the given numerical sequence and on Hahn's theorem}, Izv.Akad.Nauk, {\bf 4}
(1940), 215-228 (in Russian).


\bi{Hahn} W.Hahn, {\it \"Uber die Jacobischen Polynome und Zwei
verwandte Polynomklassen}, Math.Z. {\bf 39} (1935), 634-638.




\bi{Ismail} M.E.H.Ismail, {\it Classical and Quantum orthogonal polynomials in one variable}.
Encyclopedia of Mathematics and its Applications (No. 98), Cambridge, 2005.


\bi{Khol} A.N.Kholodov, {The umbral calculus and orthogonal polynomials}, Acta Appl.Mathem. {\bf 19} (1990), 1--54.

\bibitem{KLS} R. Koekoek, P.A. Lesky, and R.F. Swarttouw. {\it Hypergeometric orthogonal polynomials and their q-analogues}. Springer, 1-st edition, 2010.

\bi{KS} H.L.Krall and I.M.Sheffer, {\it On pairs of related orthogonal polynomial sets}, Math. Zeitschr. {\bf 86} (1965), 425--450.

\bi{Kwon} K.H. Kwon, L.L. Littlejohn, B.H. Yoo, {\it
Characterization of orthogonal polynomials satisfying differential
equations}, SIAM J. Math. Anal. {\bf 25} (1994) 976-�990.

\bi{KY} K.H.Kwon, G.J.Yoon, {\it Generalized Hahn's theorem}, J.Comp.Appl.Math. {\bf 116} (2000), 243–-
262.

\bi{Lando} S.K.Lando, {\it Lectures on Generating Functions}, AMS, 2003.

\bi{Little} Littlejohn L.L., {\it The Krall polynomials: a new class of orthogonal polynomials}, Quaestiones Math. {\bf 5} (1982),
255–-265.

\bi{Little2} Littlejohn L.L., {\it On the classification of differential equations having orthogonal polynomial solutions}, Ann.
Mat. Pura Appl. (4) {\bf 138} (1984), 35-–53.



\bi{Maroni} P.Maroni, {\it Une th\'eorie alg\'ebrique des polynomes orthogonaux. Application aux polyn\^omes orthogonaux semi-classiques}. In: Orthogonal Polynomials and their Applications. Eds. C. Brezinski et al., IMACS, Ann. Comp. Appl.
Math. 9 (Baltzer, Basel, 1991), pp. 95--130.


\bi{NSU} A.F. Nikiforov, S.K. Suslov, and V.B. Uvarov, {\em
Classical Orthogonal Polynomials of a Discrete Variable},
Springer, Berlin, 1991.


\bi{Roman} S.Roman, {\it The Theory of the Umbral Calculus. I }, J.Math.Anal.Appl. {\bf 87} (1982), 58--115.

\bi{SVZ} V.Spiridonov, L.Vinet and A.Zhedanov, {\it Spectral
transformations, self-similar reductions and orthogonal polynomials}, J.Phys.
A:  Math.  and Gen.  {\bf 30} (1997), 7621-7637.

\bi{Spi}  Spiridonov V.P., {\it Essays on the theory of elliptic hypergeometric functions}, Russ. Math. Surv. {\bf 63} (2008),
405-–472, arXiv:0805.3135.


\bi{TZ_dense} S.Tsujimoto and A.Zhedanov, {\it Elliptic Hypergeometric Laurent Biorthogonal Polynomials with a Dense Point Spectrum on the Unit Circle}, SIGMA {\bf 5} (2009), 033. arXiv:0809.2574.

\bi{VYZ} L.Vinet, O.Yermolayeva and A.Zhedanov, {\it A method to
study the Krall and q-Krall polynomials}, J.Comp.Appl.Math. {\bf
133} (2001) 647-�656.

\bibitem{VZ_missing}
L.~Vinet and A.~Zhedanov.
\newblock {A `missing' family of classical orthogonal polynomials}.
\newblock {\em J. Phys. A: Math. Theor.}, 44:085201, 2011.


\bi{WW} E.T. Whittacker, G.N. Watson, {\em A Course of Modern
Analysis}, Cambridge, 1927.

\bi{ZheR} A.Zhedanov, {\it Rational spectral transformations and orthogonal polynomials}, J.Comput.Appl.MAth. {\bf 85} (1997), 67--86.


\bi{Zhe_hyp} A.Zhedanov, {\it Abstract "hypergeometric" orthogonal polynomials}, arXiv:1401.6754.

\eb

\end{document}